# Quasi Fibonacci approximation to the low tiny fluctuations of the Li-Keiper coefficients: a numerical computation

Danilo Merlini[1], Massimo Sala[2] and Nicoletta Sala[3]

[1] *CERFIM/ISSI, Locarno, Switzerland, e-mail: merlini@cerfim.ch*
[2] *Independent Researcher*
[3] *Institute for the Complexity Studies, Rome, Italy*

**Abstract**

Using the first discrete derivatives for the expansion in z=0 of the oscillating part $\lambda_{tiny}(n) = \lambda_n^*$ of the "tiny" Li-Keiper coefficients, we analyse two series in the variable z=1-1/s ~0 for the first low values and compare them with the exact series. The numerical results suggest interesting more "sophisticated" approximations.

## 1. Introduction

In this short note subsequent to some recent works [1, 2, 3] we carry out two simple numerical experiments concerning the first few terms of the oscillating part of the Li-Keiper coefficients and compare those coefficients with the true values emerging from the exact Taylor expansion in the variable $z=1-1/s$.

## 2. The two approximations

We consider (for Re(s) ~1) the function

$$f(s) = s \cdot (s-1) \cdot d/ds \, (\log[(s-1) \cdot \zeta(s)]) \qquad (1)$$

where $\zeta(s)$ is the Zeta function of argument s ; after a change of the variable $s = 1/(1-z)$ i.e. $z= 1-1/s$, (here now $|z| \sim 0$ i.e. in the open disk), we have

$$f(1/(1-z)) = z/(1-z)^2 \, (dz/ds) \cdot d/dz \, \left( \sum_{n=1}^{\infty} (\lambda_{tiny}(n)/n) \cdot z^n \right) =$$

$$= z \cdot \left( \sum_{n=1}^{\infty} n \cdot (\lambda_{tiny}(n)/n) \cdot z^{n-1} \right) = \sum_{n=1}^{\infty} (\lambda_{tiny}(n) \cdot z^n \qquad (2)$$

Eq.(2) follows from the expansion of the log of the $\xi(s)$ function, i.e.

$$\log(\xi(s)) = \log(1/2) + \sum_{n=1}^{\infty} [(\lambda_{trend}(n) + \lambda_{tiny)}(n))/n] \cdot z^n, \qquad (3)$$

$$\log((s-1).\zeta(s)) = \log(z/(1-z))\zeta(1/(1-z)) = \sum_{n=1}^{\infty} [\lambda_{tiny}(n)/n] \cdot z^n,$$

$(z=1-1/s)$, $s=1/(1-z)$ and $d/ds = dz/ds$, $d/dz = (1-z)^2 \cdot d/dz$
let $\lambda_{tiny)} = \lambda^*_n$ ; we are interested to study the function f , i.e.

$$f(z) = d/dz \, (\log(z/(1-z))\zeta(1/(1-z))) = \sum_{n=1}^{\infty} [\lambda_{tiny}(n)] \cdot z^{n-1} \cdot (1-z)^2,$$

that is

$$f(z) = \lambda_1^* + (\lambda_2^* - 2\cdot\lambda_1^*)\cdot z + (\lambda_3^* - 2\cdot\lambda_2^* + \lambda_1^*)\cdot z^2 + (\lambda_n^* - 2\cdot\lambda_{n-1}^* + \lambda_{n-2}^*)\cdot z^3 +$$
$$+ \ldots..$$

and we analyse the sequence $\varphi(n) = (\lambda_n^* - 2\cdot\lambda_{n-1}^* + \lambda_{n-2}^*)$, $n \geq 3$,
where $\varphi(1) = \gamma = 0.577\ldots$ is the Euler-Mascheroni constant
$\varphi(1) = \lambda_1^* = \gamma = 0.577\ldots$; $\varphi(2) = \lambda_2^* - 2\cdot\lambda_1^*$ and for $n \geq 3\ldots$
$\varphi(n) = (\lambda_n^* - 2\cdot\lambda_{n-1}^* + \lambda_{n-2}^*)$, related to the constant first order discrete
derivative i.e. from the relation $\lambda_n^* - \lambda_{n-1}^* = \lambda_{n-1}^* - \lambda_{n-2}^*$.
The discrete function $\varphi(n)$ of the variable $n$ is given below up
to $n=33$ from the exact Taylor expansion around $z=1-1/s \sim 0$ ($s \sim 1$).

Notice, $f(x)$ is the function $\varphi(x=n)$, constant in the interval $(n, n+1)$.

$$\begin{aligned}
f := x \to\; & 0.5772156649 \cdot \text{Heaviside}(1 - x) \\
& - (0.1875462328) \cdot \text{Heaviside}(x - 1) \\
& \cdot \text{Heaviside}(2 - x) - (0.1358576008) \\
& \cdot \text{Heaviside}(x - 2) \cdot \text{Heaviside}(3 - x) \\
& - (0.09892062760) \cdot \text{Heaviside}(x - 3) \\
& \cdot \text{Heaviside}(4 - x) - (0.07221083529) \\
& \cdot \text{Heaviside}(x - 4) \cdot \text{Heaviside}(5 - x) \\
& - (0.05265054126) \cdot \text{Heaviside}(x - 5) \\
& \cdot \text{Heaviside}(6 - x) - (0.03813731388) \\
& \cdot \text{Heaviside}(x - 6) \cdot \text{Heaviside}(7 - x) \\
& - (0.02722760867) \cdot \text{Heaviside}(x - 7) \\
& \cdot \text{Heaviside}(8 - x) - (0.01892424162) \\
& \cdot \text{Heaviside}(x - 8) \cdot \text{Heaviside}(9 - x) \\
& - (0.01253338388) \cdot \text{Heaviside}(x - 9) \\
& \cdot \text{Heaviside}(10 - x) - (0.007568185360) \\
& \cdot \text{Heaviside}(x - 10) \cdot \text{Heaviside}(11 - x) \\
& - (0.003683619661) \cdot \text{Heaviside}(x - 11) \\
& \cdot \text{Heaviside}(12 - x) - (0.000632399682) \\
& \cdot \text{Heaviside}(x - 12) \cdot \text{Heaviside}(13 - x) \\
& + (0.0017650238) \cdot \text{Heaviside}(x - 13) \\
& \cdot \text{Heaviside}(14 - x) + (0.00364092579) \\
& \cdot \text{Heaviside}(x - 14) \cdot \text{Heaviside}(15 - x) \\
& + (0.0050942093) \cdot \text{Heaviside}(x - 15) \\
& \cdot \text{Heaviside}(16 - x) + (0.00620030462) \\
& \cdot \text{Heaviside}(x - 16) \cdot \text{Heaviside}(17 - x) \\
& + (0.00701748948) \cdot \text{Heaviside}(x - 17) \\
& \cdot \text{Heaviside}(18 - x) + (0.00759192916) \\
& \cdot \text{Heaviside}(x - 18) \cdot \text{Heaviside}(19 - x) \\
& + (0.0079604621) \cdot \text{Heaviside}(x - 19) \\
& \cdot \text{Heaviside}(20 - x) + (0.00815371201) \\
& \cdot \text{Heaviside}(x - 20) \cdot \text{Heaviside}(21 - x) \\
& + (0.00819680494) \cdot \text{Heaviside}(x - 21) \\
& \cdot \text{Heaviside}(22 - x) + (0.008111755940) \\
& \cdot \text{Heaviside}(x - 22) \cdot \text{Heaviside}(23 - x) \\
& + (0.00791186973) \cdot \text{Heaviside}(x - 23) \\
& \cdot \text{Heaviside}(24 - x) + (0.007618863680) \\
& \cdot \text{Heaviside}(x - 24) \cdot \text{Heaviside}(25 - x) \\
& + (0.00725238578) \cdot \text{Heaviside}(x - 25) \\
& \cdot \text{Heaviside}(26 - x) + (0.006808768599) \\
& \cdot \text{Heaviside}(x - 26) \cdot \text{Heaviside}(27 - x) \\
& + (0.00631704956) \cdot \text{Heaviside}(x - 27) \\
& \cdot \text{Heaviside}(28 - x) + (0.005761042767) \\
& \cdot \text{Heaviside}(x - 28) \cdot \text{Heaviside}(29 - x) \\
& + (0.005196008338) \cdot \text{Heaviside}(x - 29) \\
& \cdot \text{Heaviside}(30 - x) + (0.004590878136) \\
& \cdot \text{Heaviside}(x - 30) \cdot \text{Heaviside}(31 - x) \\
& + (0.00397031974) \cdot \text{Heaviside}(x - 31) \\
& \cdot \text{Heaviside}(32 - x) + (0.003519405964) \\
& \cdot \text{Heaviside}(x - 32) \cdot \text{Heaviside}(33 - x);
\end{aligned}$$

Below we present the plots in intervals.

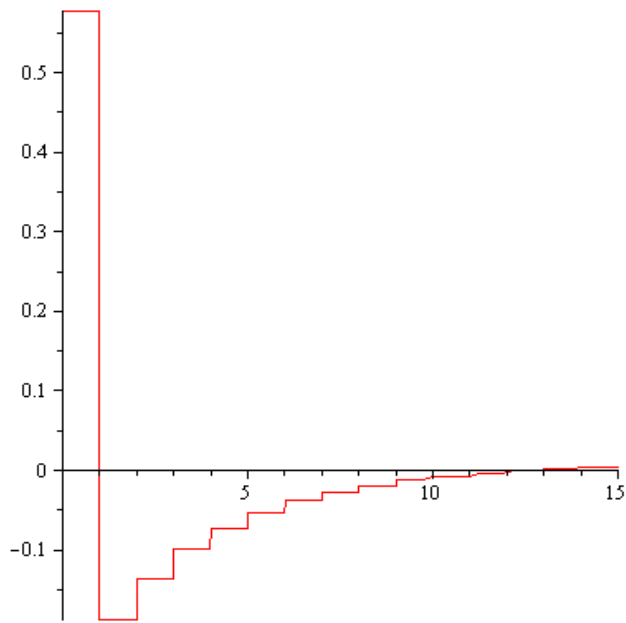

Fig. 1 The function φ(n), n=1...15.

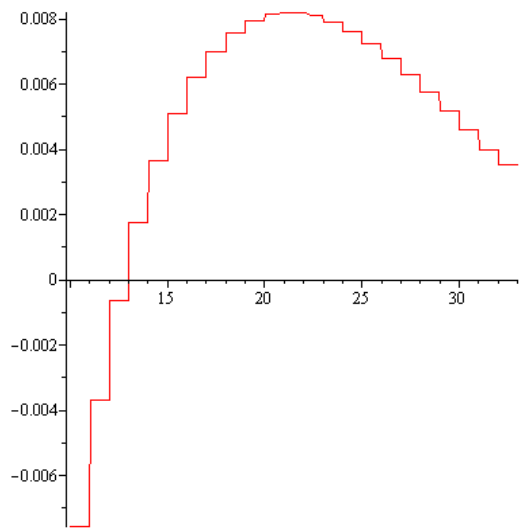

Fig. 2 The function φ(n), n=10...33.

We remark that φ(n) change sign at n=13 and then seems to decay to zero.
If this is the case, for large n we are in front of a quasi "Fibonacci" sequence

(so named in honour to the Italian mathematician Leonardo Pisano (c. 1170-1235), known as Fibonacci, derived by "filius Bonacci" ("son of Bonacci") [4]. This nickname was made up in 1838 by the Franco-Italian historian Guglielmo Libri Carucci dalla Sommaja [5]). The sequence used here is:

$\varphi(n) = (\lambda_n^* - 2\cdot\lambda_{n-1}^* + \lambda_{n-2})$, we have, with $\chi^*(n) := \lambda_{tiny}/n = \lambda_n^*/n$, and assuming $\varphi(n) \sim 0$ that:

$$\chi^*(n) = 2\cdot\chi^*(n-1) - \chi^*(n-1)$$

In fact, dividing by n

$$\chi^*(n) = \lambda_n^*/n = 2\lambda_{(n-1)}^*/(n-1+1) - \lambda_{(n-2)}^*/(n-2+2) =$$

$$\sim 2.\ \lambda_{(n-1)}^*/(n-1+1) - \lambda_{(n-2)}^*/(n-2+2)$$
$$\sim 2.\ \lambda_{(n-1)}^*/(n-1) - \lambda_{(n-2)}^*/(n-2)$$

i.e. $\quad\quad a_n \quad = \quad 2\cdot a_{n-1} \quad - \quad a_{n-2} \quad\quad (4)$

where $a_n = \lambda_{tiny}(n)/n$ (notice that, $\chi^*(n) = \lambda_{tiny}(n)/n$). Below, we give the Table of $\chi^*(n) = \lambda_{tiny}(n)/n$, which are the coefficients of the expansion of

$$\log((s-1)\cdot\zeta(s)) = \log((z/(1-z)\zeta(s)) = \sum_{n=1}^{\infty} [\lambda_{tiny}(n)/n]\cdot z^n \quad\quad (5)$$

$0.5772156649015329\, z + 0.4834425484813502\, z^2$
$+ 0.4068989760722319\, z^3$
$+ 0.3438970329678144\, z^4$
$+ 0.2916537000394335\, z^5$
$+ 0.2480497212020363\, z^6$
$+ 0.2114558343198340\, z^7$
$+ 0.1806069680149283\, z^8$
$+ 0.1545107118656992\, z^9$
$+ 0.1323803683696288\, z^{10}$
$+ 0.1135857068929158\, z^{11}$
$+ 0.09761652057825469\, z^{12}$
$+ 0.08405548593188946\, z^{13}$
$+ 0.07255781518502553\, z^{14}$
$+ 0.06283589473908967\, z^{15}$
$+ 0.05464760383629734\, z^{16}$
$+ 0.04778736544738263\, z^{17}$
$+ 0.04207923794311240\, z^{18}$
$+ 0.03737154075671643\, z^{19}$
$+ 0.03353264064005859\, z^{20}$
$+ 0.03044762189425016\, z^{21}$
$+ 0.02801563447459375\, z^{22}$
$+ 0.02614776549716923\, z^{23}$
$+ 0.02476531766445986\, z^{24}$
$+ 0.02379840619775804\, z^{25}$
$+ 0.02318480677218146\, z^{26}$
$+ 0.02286900254417814\, z^{27}$
$+ 0.02280139009413512\, z^{28}$
$+ 0.02293761303016977\, z^{29}$
$+ 0.02323799870208639\, z^{30} + O(z^{31})$

Taylor expansion of the true function Eq.(5).

It may be argued that the quasi Fibonacci sequence may furnish more values of the fluctuations (n>30) but they need then to be compared with the true values obtained from the Taylor sequence for n>30. Of course, if we start with the first three exact values of the exact series and continue, the sequence will be linear and after some $n_0$ the terms becomes negative (straight line) at around $n_0$. We now consider another approximation i.e. a second order one given by - still using- the Formulas:

$$\log((s-1)\cdot\zeta(s)) = \log(z/(1-z))\zeta(1/(1-z)) = \sum_{n=1}^{\infty} [\lambda_{tiny}(n)/n] \cdot z^n,$$

i.e.

$$(1/s)\cdot d/ds(\log((s-1)\cdot\zeta(s))) = \sum_{n=1}^{\infty} [\lambda_{tiny}(n)/n] \cdot z^{n-1}(n)\cdot(1-z)^3 =$$

$$= \sum_{n=1}^{\infty} [\lambda^*(n)] \cdot z^{n-1} \cdot (1-z)^3 \qquad (6)$$

This function differs from the above one only by the factor $(1-z)^3$ instead of the factor $(1-z)^2$: we also analyse the series where now we have, using the Taylor series above, the following expansion:

$$\lambda_1^* + (\lambda_2^* - 2\cdot\lambda_1^*)\cdot z + (\lambda_3^* - 3\cdot\lambda_2^* + \lambda_1^*) + (\lambda_4^* - 3\cdot\lambda_3^* + 3\cdot\lambda_2^* - \lambda_1^*) + ..$$

and we will analyse the functions with the help of $\varphi(n)$ where as for the first case we set

$$\varphi(n) = (\lambda_n^* - 3\cdot\lambda_{n-1}^* + 3\cdot\lambda_{n-2}^* - \lambda_{n-3}^*) \sim 0 \; ; \text{ we obtain in this way}$$
an additional sequence given by:

$$\lambda_n^* = 3\cdot\lambda_{n-1}^* - 3\cdot\lambda_{n-2}^* + \lambda_{n-3}^* \qquad n > 3. \qquad (7)$$

(here we have the second discrete derivative (applied to the $\lambda_n^*$'s)
As before we divide by n both members of the Equation and discard smaller terms; then

$$\chi^*(n) = \lambda_n^*/(n) = 3\cdot\lambda_{n-1}^*/(n-1) - 3\cdot\lambda_{n-2}^*/(n-2) + \lambda_{n-3}^*/(n-3) \qquad (8)$$

The above sequence should give values comparable to those obtained with the quasi Fibonacci sequence studied before. Notice that here the relation above involves three values of antecedents lambda's.
Below we give the Table of the numerical results of the above approximation (with 2 terms, A) and the one already derived (with 3 terms, B), but for the quantity $\lambda_{tiny}(n)/n$, up to n= 30.

We give the first 6 decimals but the calculations was carried out with 20 digits for the exact Taylor series of $\log[(z/(1-z)\cdot\zeta(1/(1-z)]$ in z=0, (C). It should be noticed that we have introduced the true values in the expression for the $\chi_n^* = \lambda_{tiny}(n)/n \cdot (\lambda_{tiny}(1)/1 = \gamma = 0.577215...)$.

| n | $\lambda_{tiny}(n)/n$ (A) | $\lambda_{tiny}(n)/n$ (C) | $\lambda_{tiny}(n)/n$ (B) |
|---|---|---|---|
| 2 | - | 0.483442 | |
| 3 | 0.452184 | 0.406898 | - |
| 4 | 0.368627 | 0.343897 | 0.334662 |
| 5 | 0.306095 | 0.291653 | 0.286311 |
| 6 | 0.256824 | 0.248049 | 0.244789 |
| 7 | 0.216904 | 0.211455 | 0.209382 |
| 8 | 0.184010 | 0.180606 | 0.179243 |
| 9 | 0.156613 | 0.154510 | 0.153588 |
| 10 | 0.133633 | 0.132380 | 0.131741 |
| 11 | 0.114273 | 0.113585 | 0.113134 |
| 12 | 0.097923 | 0.097616 | 0.097292 |
| 13 | 0.084204 | 0.084055 | 0.083820 |
| 14 | 0.072431 | 0.072557 | 0.072386 |
| 15 | 0.062593 | 0.062835 | 0.062710 |
| 16 | 0.054329 | 0.054647 | 0.054556 |
| 17 | 0.047422 | 0.047787 | 0.047722 |
| 18 | 0.041689 | 0.042079 | 0.042033 |
| 19 | 0.036971 | 0.037371 | 0.037341 |
| 20 | 0.033134 | 0.033532 | 0.033514 |
| 21 | 0.030059 | 0.030447 | 0.030438 |
| 22 | 0.027643 | 0.028015 | 0.028013 |
| 23 | 0.025795 | 0.026147 | 0.026151 |
| 24 | 0.024435 | 0.024765 | 0.024773 |
| 25 | 0.023493 | 0.023798 | 0.023810 |
| 26 | 0.022905 | 0.023184 | 0.023199 |
| 27 | 0.022616 | 0.022869 | 0.022885 |
| 28 | 0.022575 | 0.022801 | 0.022829 |
| 29 | 0.022738 | 0.022937 | 0.022956 |
| 30 | 0.023064 | 0.023237 | 0.023257 |
| | | | 0.023686 |

Table

(A) The quasi Fibonacci computed with the true values of the Taylor expansion; (C) The exact values by means of the Taylor expansion (B) the second approximation with three antecedent terms, to 6 digits.

Below we construct the plots of the above three discrete functions by intervals.

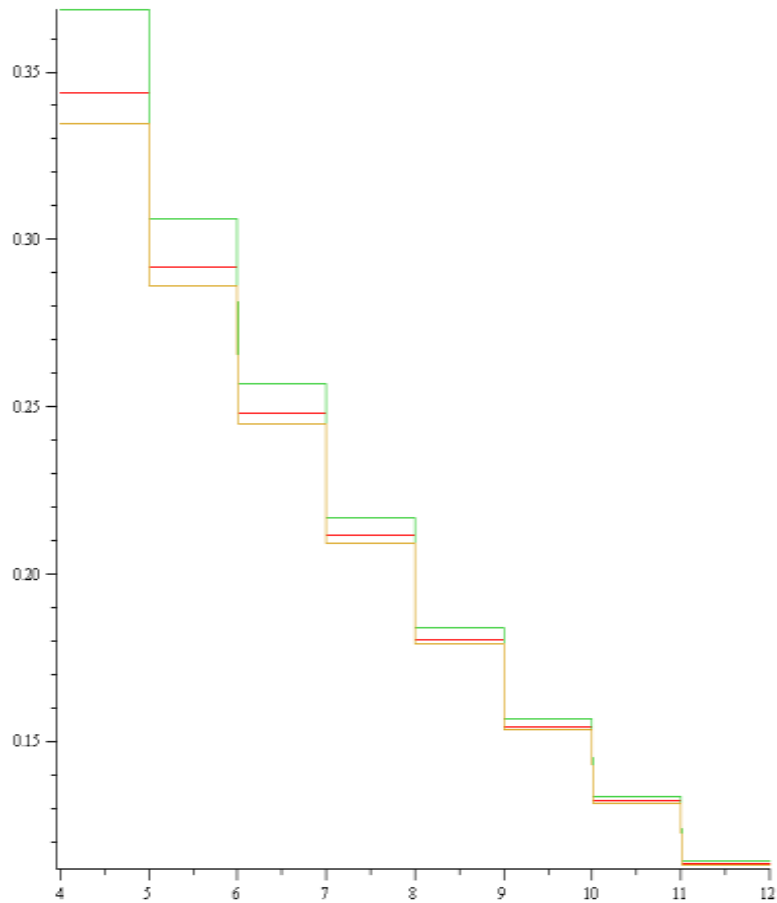

Fig. 3 Plots of the three functions in the range n= [4 .. 13]. In reed the true function, (A) in green, (B) in maroon and (C) in red.

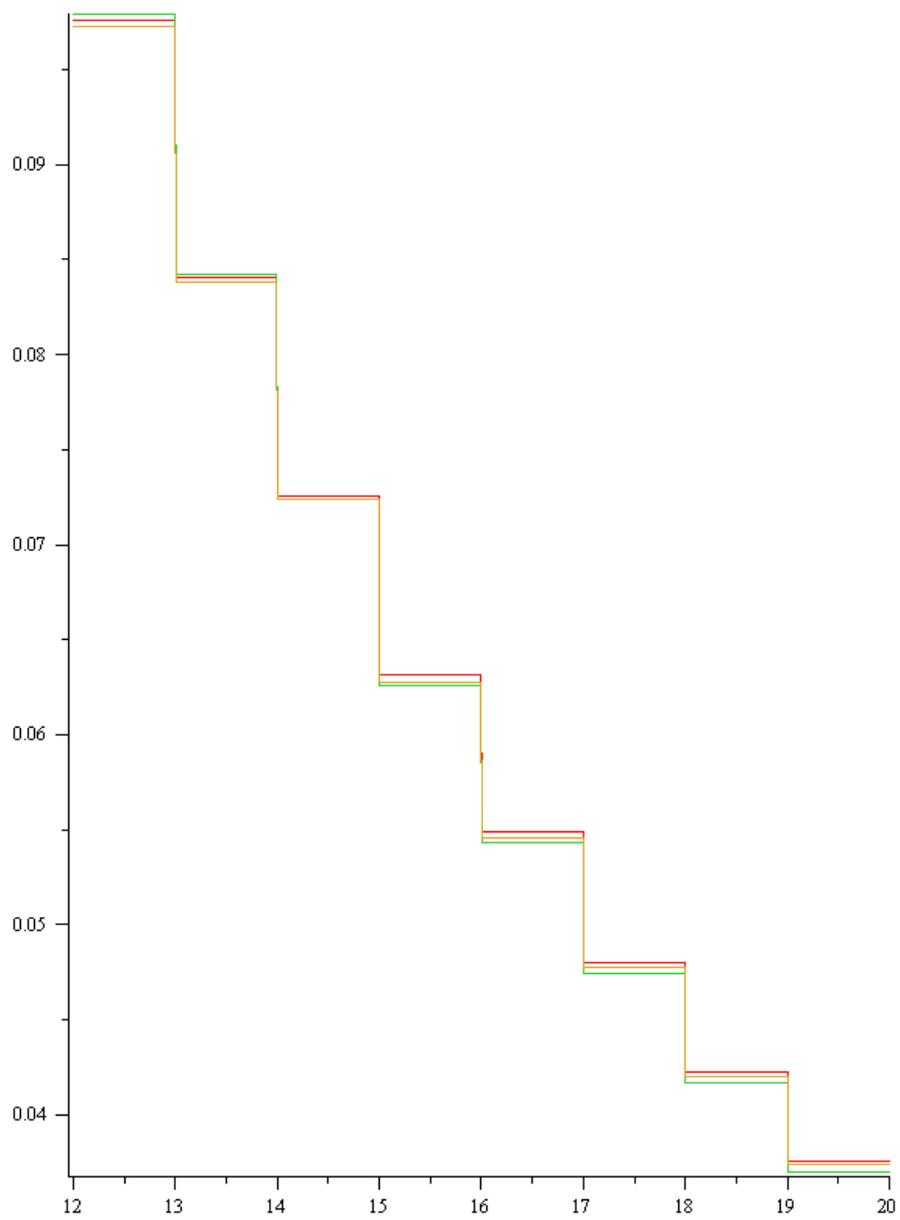

Fig. 4   The three functions in the intervall  n= [ 12, 20 ].
(A) in green, (B) in maroon and (C) in red.

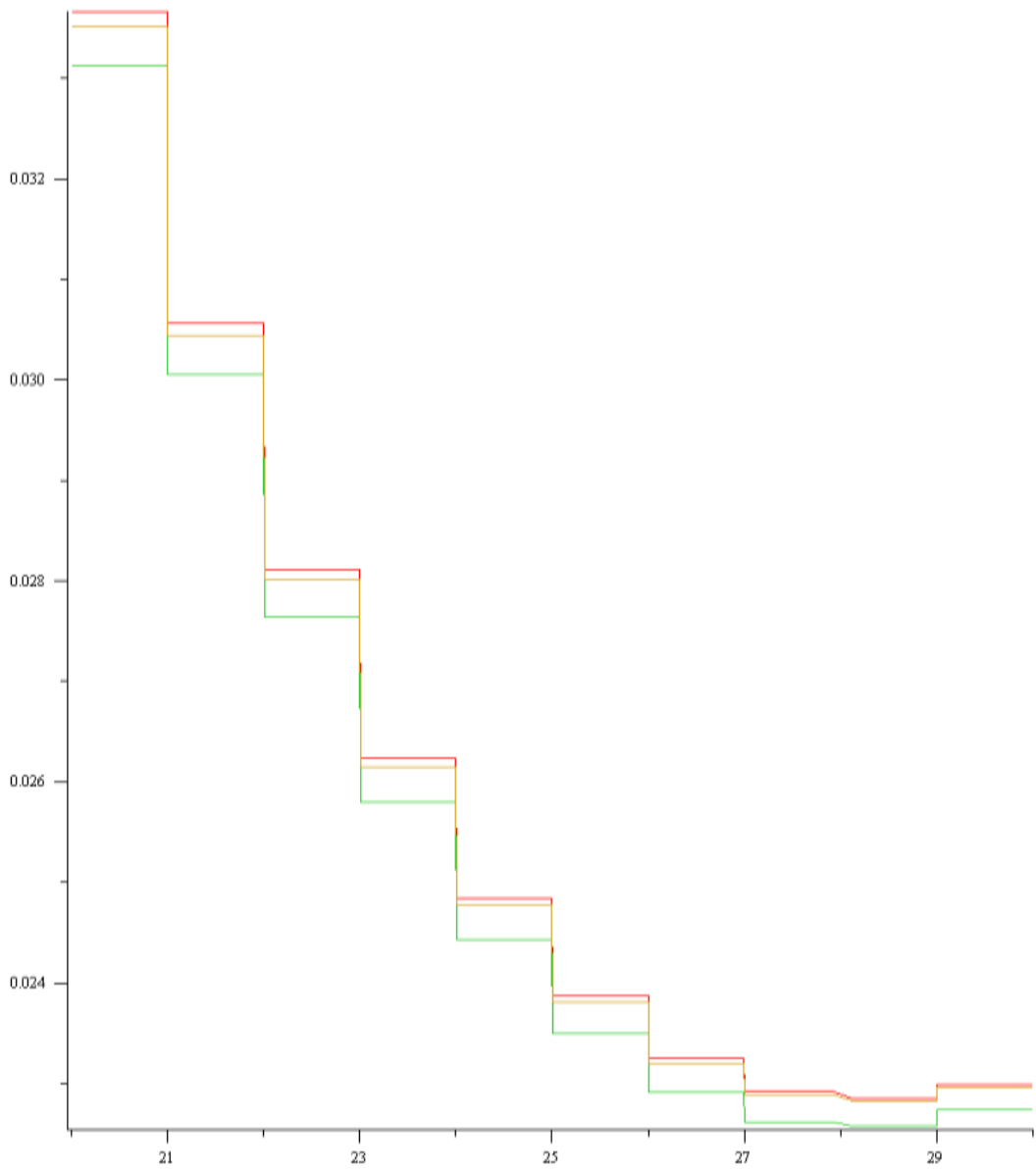

Fig. 5 The three functions in the range n = [ 20,30]
(A) in green, (B) in red and (C) in maroon.

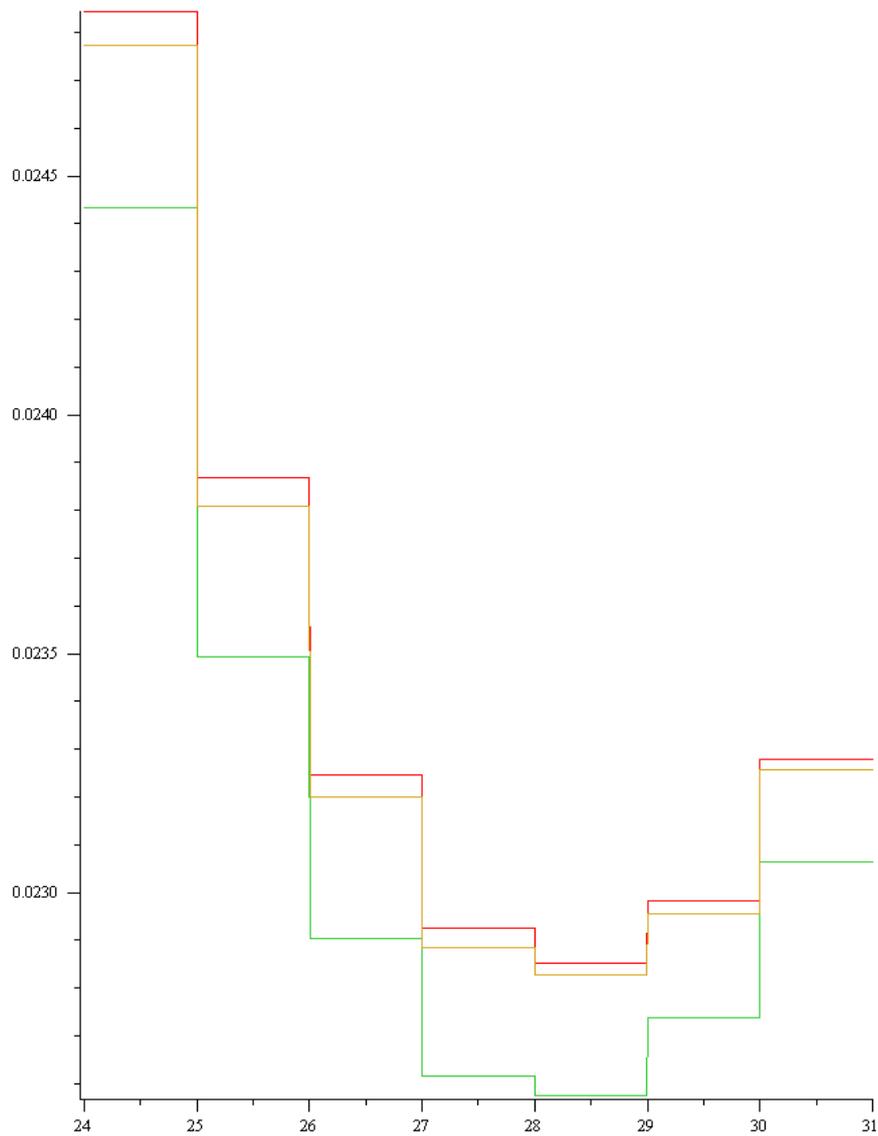

Fig.6
The plots of the three functions, (A) in green, (B) in red and (C) in maroon. The same as the Figure above but on another scale. The functions are slowly increasing around n=20, a manifestation of a possible emergence of small oscillations for bigger values of n.

Of course, if we were to use only the first two exact value $\lambda_{tiny}(n)/n$, n=1 and n=2 the corresponding straight line i.e. with  0.483442 and 0.406898 from column (C) above, has the equation:

$$\lambda_{tiny}(n)/n = \lambda_{tiny}(1)/(1) + (0.406898-0.483442)\cdot(n-1)$$

is vanishing for  n= 1+ γ/0.076544 = 1+ 0.577215/0.076544= ~ 8

If, on the other hand, we use the exact result for n=20 and n=21 and of column (C) above, then $\lambda_{tiny}(n)/n$ will vanish around n=30 and so on.

## 3. Concluding Remark

In this short note, as a preliminary analysis of a more complete work where we will introduce and analyse a more complex approximation around a special K function [6] we have computed two series related to the exact Taylor expansion for the reduced "tiny" oscillations up to n=30 , i.e. $\lambda_{tiny}(n)/n$.